\documentclass[11pt,amscd,amssymb,verbatim]{amsart}
\usepackage{amssymb,amsmath,amsthm,mathrsfs,graphicx,xcolor,tikz}
\usepackage{amsmath,amssymb,amsfonts,amsthm,graphicx}
\usepackage[mathscr]{euscript}
\usepackage{epsfig}
\usepackage{subfigure}
\usepackage{graphicx}
\usepackage{amsbsy}
\usepackage{cite}
\usepackage{amsmath}
\makeatletter
\renewenvironment{proof}[1][\proofname]{\par
\pushQED{\qed}%
\normalfont \topsep6\p@\@plus6\p@\relax
\trivlist
\item\relax
{\bfseries
#1\@addpunct{.}}\\
}{%
\popQED\endtrivlist\@endpefalse
}
\makeatother
\usepackage{color}
\usepackage{natbib}
\usepackage{nameref}
\theoremstyle{definition}
 \usepackage[colorlinks,urlcolor=blue]{hyperref}
    \usepackage{url}
\newtheorem{exmp}{Example}[section]
\usepackage{xspace}

\numberwithin{equation}{section}
\newtheorem{theorem}{Theorem}

\newtheorem{remark}{Remark}

\begin{document}

\title[Varextropy  of doubly truncated random variable] {Varextropy  of doubly truncated random variable \vspace{1cm}}
 \maketitle
 \begin{center}
\author{ R. Zamini\footnote{Corresponding author: {\tt rahelehzamini1361@yahoo.com}} $^1$,  S. Ghafouri$^2$  and F. Goodarzi$^3$\\
	\small $^1$  Department  of Mathematics, Faculty of Mathematical Sciences and  Computer, Kharazmi University, Thehran, Iran. \\
	\small $^2$   Department of Mathematics, Faculty of Science, Arak University, Arak, 384817758,\ Iran.\\
\small $^3$   Department of Statistics, Faculty of Mathematical Sciences, University of Kashan, Kashan,\ Iran.\\}
\end{center}
\date{}

\vspace{0.09 in}

\begin{abstract}
Recently, there has been growing attention to study uncertainty measures for doubly truncated random variables.
In this paper, the concept of varextropy for doubly truncated random variables is introduced. The changes of this measure under linear transforms are investigated.
Also, lower and upper bounds for the interval varextropy under some conditions are obtained and some other properties for this measure are proved.
\noindent
In the following, three estimators for the interval varextropy are proposed. A simulation study has been carried out to investigate the behavior of estimators. At the end, an application for these estimators is proposed.
\end{abstract}

\vspace{0.09 in}

\section{\bf Introduction and Preliminaries}\label{section1}
Uncertainty or randomness is associated with any event or phenomenon. A basic measure of uncertainty for a random variable was given by \cite{Shannon1948}.
Let $X$ be a non-negative random variable with the density function $f$, the Shannon entropy is defined as 
\begin{equation}\label{eq1}
H(x)=-\int_{0}^{\infty}  f(x)\log(f(x))\,d{x}.\\
\end{equation}
\cite{Song2001} introduced the varentropy function. This   measure  is the expectation of the squared deviation of the information content $-\log(f(x))$ from entropy:
\begin{equation}\label{eq2}
VJ(X)=Var(-\log(f(X))=\int_{0}^{\infty}  f(x)(\log(f(x)))^2\,d{x}-\bigg[\int_{0}^{\infty}  f(x) \log(f(x))\bigg]^2\,d{x}.\\
\end{equation}
$VJ(X)$ measures how much variability we may expect in the information  content of $X$.  

\cite{DiCrescenzoandPaolillo2021} defined the varentropy for the residual lifetime $X_{t}=(X-t|X>t)$, called the residual varentropy at time $t$ and defined as 
\begin{equation}\label{eq3}
VJ(X_t)=\int_{t}^{\infty}  \frac{f(x)}{\overline{F}(t)}\bigg(\log\big(\frac{f(x)}{\overline{F}(t)}\big)\bigg)^2\,d{x}-\bigg[\int_{t}^{\infty}  \frac{f(x)}{\overline{F}(t)} \log\big(\frac{f(x)}{\overline{F}(t)}\big)\,d{x}\bigg]^2,\\
\end{equation}
where $\overline{F}(t)=1-F(t)$.

\cite{BuonoandLongobardi2020} studied the dual measure of residual varentropy for  the past lifetime $X_{(t)}=(t-X|X \leq t)$, called past varentropy and defined as follows:
\begin{equation}\label{eq4}
VJ(X_{(t)})=\int_{0}^{t}  \frac{f(x)}{F(t)}\bigg(\log\big(\frac{f(x)}{F(t)}\big)\bigg)^2\,d{x}-\bigg[\int_{0}^{t}  \frac{f(x)}{F(t)} \log\big(\frac{f(x)}{F(t)}\big)\,d{x}\bigg]^2.\\
\end{equation}
The dual of entropy is called the extropy and is defined as:
\begin{equation}\label{eq5}
J(X)=-\frac{1}{2}\int_{0}^{\infty} (f(x))^2\,d{x}.\\
\end{equation}
A detailed motivation for $J(X)$ is explained in \cite{Ladetal2015}. Among other authors who have investigated the properties of the extropy and its applications, we can mention \cite{Qui2017}, \cite{RaqabandQiu2019}, \cite{Qiuetal2019}, \cite{MartinasandFrankowicz2000}, \cite{Becerraetal2018} and \cite{Balakrishnanetal2022}.

\cite{QiuandIia2018} defined the concept of measuring uncertainty in residual lifetimes called the residual extropy:
\begin{equation}\label{eq6}
J(X_t)=-\frac{1}{2({\overline{F}(t)})^2}\int_{t}^{\infty} (f(x))^2\,d{x}.\\
\end{equation}
\cite{Krishnanetal2020} and \cite{KamariandBuono2021} introduced the notion of the past extropy as 
\begin{equation}\label{eq7}
J(X_{(t)})=-\frac{1}{2({F(t)})^2}\int_{0}^{t} (f(x))^2\,d{x},\\
\end{equation}
which measures the uncertainty in the time $X_{(t)}$. An alternative measure to \eqref{eq2} is introduced by \cite{Vaselabadietal2021}. This measure is known as varextropy and defined as 
\begin{equation}\label{eq8}
V(X)=Var(-\frac{1}{2}f(X))=\frac{1}{4}\int   (f(x))^3 \,d{x}-\bigg(\frac{1}{2}\int (f(x))^2\,d{x}\bigg)^2.\\
\end{equation}
This measure indicates how the information content is scattered around the extropy. \cite{NoughabiandNoughabi2024} proposed some estimators for the varextropy of a continuous random variable.

\cite{Zaminietal2023} proposed non-parametric estimators for varextropy under length-biased sampling. They investigated some asymptotic properties of the estimators such as a.s. consistency and asymptotic normality. \cite{Goodarzi2022} showed that, if extropy and varextropy series and parallel systems are equal then distribution function $F(x)$ is a symmetric function.

For measuring the dispersion of the residual and past lifetimes, \cite{Vaselabadietal2021} introduced $V(X_t)$ and $V(X_{[t]})$ as follows:
\begin{equation}\label{eq9}
V(X_t)=\frac{1}{4(\overline{F}(t))^3}\int_{t}^{\infty}  (f(x))^3 \,d{x}-\bigg(\frac{1}{2(\overline{F}(t))^2}\int_{t}^{\infty} (f(x))^2\,d{x}\bigg)^2,\\
\end{equation}
\begin{equation}\label{eq10}
V(X_{[t]})=\frac{1}{4(F(t))^3}\int_{0}^{t}  (f(x))^3 \,d{x}-\bigg(\frac{1}{2(F(t))^2}\int_{0}^{t} (f(x))^2\,d{x}\bigg)^2.\\
\end{equation}
Doubly truncated data play an important role in the many fields such as reliability theory, economy and astronomy. Recently, uncertainty measures for doubly truncated random variable are studied.
Let $X_{t_1,t_2}=(X|t_1<X<t_2)$ be the lifetime of an item falls in the interval $(t_1,t_2)$. \cite{SharmaandKundu2023}  introduced the concept of varentropy for $X_{t_1,t_2}$. They  applied this concept to some known distributions. They also studied  the behavior of the interval varentropy in terms of varied transformations. The interval varentropy is introduced as 
\begin{eqnarray}\label{eq11}\nonumber
VJ(X_{t_1,t_2})&=&\int_{t_1}^{t_2}  \frac{f(x)}{F(t_2)-F(t_1)}\bigg(\log\big(\frac{f(x)}{F(t_2)-F(t_1)}\big)\bigg)^2\,d{x}\\ \nonumber &&-\bigg[\int_{t_1}^{t_2}  \frac{f(x)}{F(t_2)-F(t_1)} \log\big(\frac{f(x)}{F(t_2)-F(t_1)}\big)\,d{x}\bigg]^2.\\
\end{eqnarray}
\eqref{eq11} takes the form of the residual varentropy \eqref{eq3} when  $ t_2 \rightarrow \infty$ and of the past varentropy  \eqref{eq4} when $ t_1 \rightarrow 0$.

\cite{Buonoetal2023} proposed the dynamic version of extropy for doubly truncated random variables as: 
\begin{equation}\label{eq12}
IJ(t_1,t_2)=J(X_{t_1,t_2})=-\frac{1}{2(F(t_2)-F(t_1))^2}\int_{t_1}^{t_2} (f(x))^2\,d{x}.\\
\end{equation}
They evaluated some characterizations of $J(X_{t_1,t_2})$ and presented some bounds for it. They also evaluated this measure under the effect of linear transformations.

In analogy with the interval varentropy and interval extropy, we introduce the concept of interval varextropy for doubly truncated random variables. This measure indicates how the information content is scattered around the interval extropy for a doubly truncated random variable.
Similar to varextropy, interval varextropy is used to determine which interval extropy is more appropriate for measuring uncertainty when the interval extropies of two doubly truncated variables are equal.
On the other hand, the interval varextropy is free of the model parameters in some distributions and therefore is more flexible than the interval varentropy.

In the rest of the paper, we introduce the interval varextropy $IV(t_1, t_2)$ for doubly truncated random variable.
In the following, several examples are shown to illustrate the interval extropy and the interval varextropy for some reliability distributions. Also, we investigate some properties
 of the $IV(t_1,t_2)$ and  give some bounds for the interval varextropy. In Section \ref{section3}, we propose some estimators for the interval varextropy and investigate their some asymptotic properties.
In  Section \ref{section4}, we investigate the behavior of proposed estimators by a Monte Carlo study. We give an application for proposed estimators in Section \ref{section5}. We illustrate the use  of proposed estimators in the real case  in Section \ref{section6}.
\section{\bf Interval Varextropy Measure}\label{section2}  
Let $X_{t_1,t_2}=(X|t_1<X<t_2)$ be the lifetime of a component in a system. That is, the lifetime of the component falls in the interval $(t_1,t_2)$. 
In the present paper, an alternative measure analogous to the interval varentropy for $X_{t_1,t_2}$ is proposed where indicates how the information content is scattered around the doubly truncated extropy defined in \eqref{eq12}.

Let $D=\{(u,v)\in  \Re _{+}^{2}: F(u)<F(v)\},$ the interval varextropy for $(t_1, t_2) \epsilon D$ is defined as:
\begin{eqnarray}\label{Khas}
IV(t_1,t_2)&=&var(-\frac{1}{2}f(X_{t_1,t_2}))=\frac{1}{4}\frac{1}{(F(t_2)-F(t_1))^3}\int_ {t_1}^{t_2} (f(x))^3\,d{x}\\
  \nonumber &&-\frac{1}{4}\frac{1}{(F(t_2)-F(t_1))^4}\big(\int_ {t_1}^{t_2} f^2(x)\,d{x}\big)^2,
\end{eqnarray}
which is an extension of varextropy.
\begin{remark}\label{R1}
It is easy to see that $IV(0,t_2)=V(X_{[t_2]}),\ IV(t_1,\infty)=V(X_{t_1})$ and $IV(0,\infty)=V(X)$ are  the past varextropy, residual varextropy and varextropy, respectively.
\end{remark}
In the following, a few examples are given to illustrate the interval extropy and the interval varextropy for significant distributions used in reliability and  survival analysis.
\begin{exmp}\label{Exam1}
Let $X$ be a random variable having exponential distribution, that is, $X \sim Exp(\lambda),\ \lambda>0$. For $0<t_1<t_2,$ we have $IJ(t_1,t_2)=-\frac{\lambda}{4}\frac{e^{-\lambda t_1}+e^{-\lambda t_2}}{e^{-\lambda t_1}-e^{-\lambda t_2}}$, and   $IV(t_1,t_2)=\frac{\lambda ^{2}}{48} \frac{\big(e^{-\lambda t_1}+e^{-\lambda t_2}\big)^2}{\big(e^{-\lambda t_1}-e^{-\lambda t_2}\big)^2}$.
 \end{exmp}
\begin{exmp}\label{Exam2}
Suppose that  $X$ have the Pareto-I distribution with the survival function $\overline{F}(x)=(\frac{a}{x})^b;\    0<a<x,\  b>0$.  We have, for $a<t_1<t_2$; 
$$IJ(t_1,t_2)=\frac{b^2}{2(2b+1)(t_1t_2)}\frac{(t_1^{2b+1}-t_2^{2b+1})}{(t_2^{b}-t_1^{b})^2},$$
$$IV(t_1,t_2)=\frac{b^3}{4{(t_2^b-t_1^b)}^3(t_1^2t_2^2)}\bigg(\frac{\big(t_2^{3b+2}-t_1^{3b+2}\big)}{(3b+2)}-\frac{b}{(2b+1)^2}\frac{\big(t_1^{2b+1}-t_2^{2b+1}\big)^2}{(t_2^b-t_1^b)}\bigg).$$
  \end{exmp}
  \begin{exmp}\label{Exam3}
 Let $X$ follow power distribution having $F(x)=(\frac{x}{a})^b, \ 0<x<a, \ b>0$. Then, for $0<t_1<t_2<a$, we have 
 $$IJ(t_1,t_2)=\frac{-b^2}{2(2b-1)(t_2^b-t_1^b)^2}\big(t_2^{2b-1}-t_1^{2b-1}\big),$$
$$IV(t_1,t_2)=\frac{b^3}{4{(t_2^b-t_1^b)}^3}\bigg[\frac{\big(t_2^{3b-2}-t_1^{3b-2}\big)}{(3b-2)}-\frac{b}{(2b-1)^2}\frac{\big(t_2^{2b-1}-t_1^{2b-1}\big)^2}{(t_2^b-t_1^b)}\bigg].$$
   \end{exmp}
  \begin{exmp}\label{Exam4}
Suppose that  $X$  is a random variable with distribution function  $F(x)=x^2, \  x \in (0,1)$. Then, for $0<t_1<t_2<1$, 
 $$IJ(t_1,t_2)=\frac{-2\big(t_2^2+t_1t_2+t_1^2\big)}{3(t_2^2-t_1^2)(t_2+t_1)},$$
$$IV(t_1,t_2)=\frac{1}{(t_2^2-t_1^2)^2}\bigg[\frac{1}{2}(t_2^2+t_1^2)-\frac{4}{9}\frac{\big(t_2^2+t_1t_2+t_1^2\big)^2}{\big(t_2+t_1\big)^2}\bigg].$$
   \end{exmp}
\begin{exmp}\label{Exam5}
Suppose that  $X$  is a random variable with density and distribution function 
$$F(x)=\begin {cases}
\exp(-\frac{1}{2}-\frac{1}{x}) &   if\  0 \leq x\leq 1 \\  \exp(-2+\frac{x^2}{2}) &  if\ 1 \leq x \leq 2, 
\end{cases} $$
$$f(x)=\begin {cases}
\frac{1}{x^2}\exp(-\frac{1}{2}-\frac{1}{x}) &    if\ x \in (0,1] \\ x\exp(-2+\frac{x^2}{2}) &  if\  x \in [1,2). 
\end{cases} $$
   \end{exmp}
   \vspace{5mm}
\begin{figure}[!ht]
\centering
\subfigure[{plot of interval varextropty $IV(0.5,t_2)$ against $t_2 \in (1,1.8)$}]{\label{mylabel1}
\includegraphics*[width=0.48\textwidth]{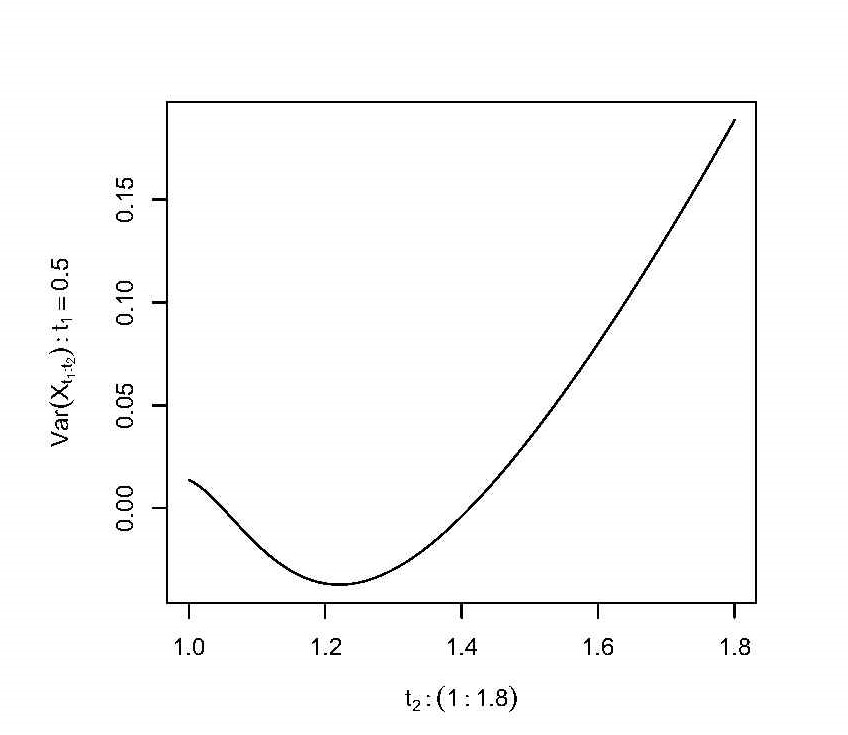}}
\hspace{0.9mm}
\subfigure[{plot of interval varextropty $IV(t_1,1.8)$ against $t_1 \in (0,1)$}]{\label{mylabel2}
\includegraphics*[width=0.48\textwidth]{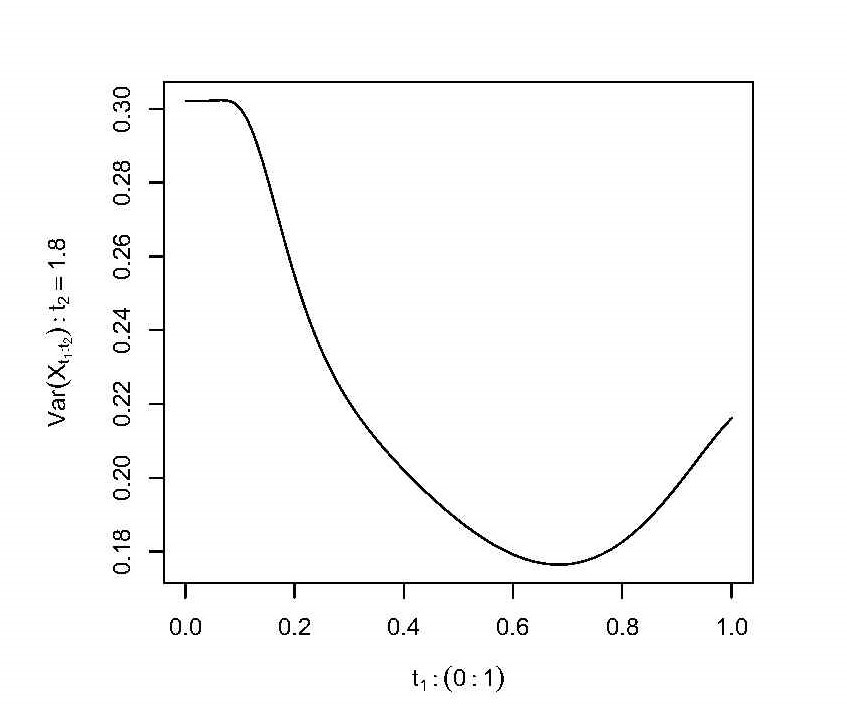}}
\caption{Graphical representation of $IV(t_1,t_2)$ in Example \ref{Exam5}}
\label{mylabel}
\end{figure}

Since the expression of the interval varextropy is not given in terms of elementry functions, we plot the interval varextropy as a function of $t_1$ for fixed $t_2$ (see Figure \ref{mylabel2}) and a function of $t_2$ for fixed $t_1$ (see Figure \ref{mylabel1}).  Figure \ref{mylabel} shows that the interval varextropy for this distribution is not monotone in $t_2$ for some fixed $t_1$ and not monotone in $t_1$ for some fixed $t_2$.
\subsection{\bf Some Properties for Interval Varextropy}\label{2.1}
In this subsection, we investigate some properties of the interval varextropy.

\noindent In the next theorem, we look for interval varextropy under a linear transformation.
\begin{theorem}\label{Theo1}
Suppose that $X$ is a non-negative and absolutely continuous random variable. Let $Y=aX+b$, with $a>0, \  b>0$. Then, their interval varextropies are related by $$IV_{Y}(t_1,t_2)=\frac{1}{a^2}IV_{X}(\frac{t_2-b}{a}, \frac{t_1-b}{a}).$$
\end{theorem}

 \begin{proof}
 The proof is similar as that of Proposition 3 of \cite{Buonoetal2023}.
 \end{proof}

\cite{NavarroandRuiz1996} defined the Generalized Failure Rate (GFR) functions of random variable $X$ in $t_1$ and $t_2$ as 
$$h_i(t_1, t_2)=\frac{f(t_i)}{(F(t_2)-F(t_1))},\ i=1,2,\ (F(t_2)-F(t_1)>0),$$
where $f$ and $F$ are density and distribution function of random variable $X$, respectively.

In the following theorem, we show that under some conditions on $IV(t_1, t_2)$, $IJ(t_1, t_2)$ can be obtained in terms of GFR functions.

\begin{theorem}\label{Theo2}
Let the interval varextropy $IV(t_1,t_2)$ be constant, that is, $IV(t_1,t_2)=v\geq 0,$ for all $(t_1,t_2) \in D$. Then, $$IJ(t_1,t_2)=-\frac{1}{4}\big(h_1(t_1,t_2)+h_2(t_1,t_2)\big)$$.
\end{theorem}

\begin{proof}
Since  $IV(t_1,t_2)=v$, is constant, we have 
\begin{eqnarray*}
\frac{\partial IV(t_1,t_2)}{\partial t_1}&=&3h_1(t_1,t_2)v-h_1(IJ(t_1,t_2))^2\\ &&-\frac{1}{4}(h_1(t_1,t_2))^3-(h_1(t_1,t_2))^2 IJ(t_1,t_2)=0,
\end{eqnarray*}
\begin{eqnarray*}
\frac{\partial IV(t_1,t_2)}{\partial t_2}&=&-3h_2(t_1,t_2)v+h_2(t_1,t_2)(IJ(t_1,t_2))^2\\ &&+(h_2(t_1,t_2))^2 IJ(t_1,t_2)+\frac{1}{4}(h_2(t_1,t_2))^3=0.
\end{eqnarray*}
By solving these two equations simultaneously, the desired result is easily obtained. 
\end{proof}

\begin{theorem}\label{Theo3}
Let  $X$ be a random variable with support $(0,\infty)$, differentiable and strictly positive density function $f$ and distribution function $F$. If $IV(t_1,t_2)=v \ge 0$  for all $(t_1,t_2) \in D$, then $X$ is exponentially distributed.
\end{theorem}

\begin{proof}
The proof of this theorem is obtained according to Theorem 2 and Theorem 4 of \cite{Buonoetal2023}.
\end{proof}

In the following, we give some bounds for the interval varextropy. Suppose that $m(t_1,t_2)$ and $\sigma ^2(t_1, t_2)$ express the mean and variance of $X$ in the interval $(t_1,t_2)$ which are defined as 
\begin{equation}\label{1}
m(t_1,t_2)=E(X|t_1<X<t_2)=\int_{t_1}^{t_2} \frac{xf(x)}{F(t_2)-F(t_1)}\,d{x},\\
\end{equation}
and
\begin{equation}\label{2}
\sigma^2(t_1,t_2)=Var(X|t_1<X<t_2)=\int_{t_1}^{t_2} \frac{x^2f(x)}{F(t_2)-F(t_1)}-(m(t_1,t_2))^2\,d{x}.\\
\end{equation}
A lower bound for $IV(t_1,t_2)$ in terms of $m(t_1,t_2)$ and $\sigma^2(t_1,t_2)$ is given in the following theorem.
\begin{theorem}\label{Theo4}
Let  $f_{t_1,t_2}(x)=\frac{f(x)}{F(t_2)-f(t_1)}$ be  the density function of random variable $X_{t_1,t_2}=(X|t_1<X<t_2)$. We have  
\begin{equation}\label{3}
IV(t_1,t_2) \geq \frac{1}{4} \sigma^2 (t_1,t_2) \bigg( E_{t_1,t_2}\bigg[( W'_{t_1,t_2}(X_{t_1,t_2})f_{t_1,t_2}(X_{t_1,t_2})\bigg]\bigg)^2,
\end{equation}
\end{theorem}
where
\begin{equation}\label{4}
 \sigma^2 ({t_1,t_2}) W_{t_1,t_2}(x)f_{t_1,t_2}(x)=\int_{0}^{x}(m(t_1,t_2)-z)f_{t_1,t_2}(z)\,d{z}.
\end{equation}

\begin{proof}
From \cite{CacoullosandPapathanasiou1989}, we can write 
\begin{equation}\label{5}
 IV(t_1,t_2)= \frac{1}{4} Var\big(f_{t_1,t_2}(X_{t_1,t_2})\big) \geq \frac{1}{4} \sigma^2 (t_1,t_2)\bigg( E_{t_1,t_2}\big(W_{t_1,t_2}(X_{t_1,t_2})\big)f'_{t_1,t_2}(X_{t_1,t_2})\bigg)^2,
\end{equation}
where $W_{t_1,t_2}(x)$ is defined as
\begin{equation}\label{6}
\sigma^2(t_1,t_2)W_{t_1,t_2}(x)f_{t_1,t_2}(x)=\int_{0}^{x}(m(t_1,t_2)-z)f_{t_1,t_2}(z)\,d{z}.\\
\end{equation}
On differentiating both sides of \eqref{6}, we can see that:
$$W_{t_1,t_2}(x)f'_{t_1,t_2}(x)=\frac{(m(t_1,t_2)-x)f_{t_1,t_2}(x)}{\sigma^2 (t_1,t_2)}-W'_{t_1,t_2}(x)f_{t_1,t_2}(x).$$
Therefore,
$$IV(t_1,t_2) \geq \frac{1}{4} \sigma^2 (t_1,t_2) \bigg( E_{t_1,t_2}\bigg[(W'_{t_1,t_2}(X_{t_1,t_2})f_{t_1,t_2}(X_{t_1,t_2})\bigg]\bigg)^2.$$
\end{proof}

Many density functions belong to the exponential family, such as the normal density function, gamma, exponential, Chi-squared, Beta, Dirichlet, etc. In the following theorem, we find an upper bound for varextropy of a doubly truncated random variable if its density belong to the exponential family.
\begin{theorem}\label{Theo5}
Let  $X$ be  a random variable such that its density function belongs to the exponential family that is 
$$f(x)=\exp\big(\eta(\theta)\tau(x)+Q(\theta)+R(x)\big),\     x \in \Re$$ 
where $\eta(\theta)$, $\tau(x)$, $Q(\theta)$ and $R(x)$ are known functions, $\eta(\theta)\neq 0$ and $R(x)$ is non-negative. Under conditions $A(t_1,t_2)=\int_{t_1}^{t_2} \tau(x)\,d{x}<\infty$ and $B(t_1,t_2)=\int_{t_1}^{t_2} R(x)\,d{x}<\infty$.

We can see that
$$IV(t_1,t_2) \leq \frac{1}{4(F(t_2)-F(t_1))^3}-\frac{1}{4(F(t_2)-F(t_1))^4}\bigg[2\eta(\theta)A(t_1,t_2)+2Q(\theta)(t_2-t_1)+2B(t_1,t_2)+(t_2-t_1)\bigg]$$
 \end{theorem}
 
 \begin{proof}
 Using the inequality $1+x< e^x$, $-\infty<x< \infty$, we can write:
 \begin{eqnarray}\nonumber
 IV(t_1,t_2)&=& \frac{1}{4(F(t_2)-F(t_1))^3}\int_{t_1}^{t_2}\exp\big(3\eta(\theta)\tau(x)+3Q(\theta)+3R(x)\big)\,d{x}\\ \nonumber &&-\frac{1}{4(F(t_2)-F(t_1))^4}\bigg(\int_{t_1}^{t_2}\exp\big(2\eta(\theta)\tau(x)+2Q(\theta)+2R(x)\big)\,d{x}\bigg)^2\\ \nonumber && \leq \frac{G(t_2)-G(t_1)}{4(F(t_2)-F(t_1))^3}-\frac{\bigg(\int_{t_1}^{t_2}\big(2\eta(\theta)\tau(x)+2Q(\theta)+2R(x)+1\big)\,d{x}\bigg)^2}{4(F(t_2)-F(t_1))^4}\\ \nonumber && \leq \frac{1}{4(F(t_2)-F(t_1))^3}-\frac{\bigg[2\eta(\theta)A(t_1,t_2)+2Q(\theta)(t_2-t_1)+2B(t_1,t_2)+(t_2-t_1)\bigg]}{4(F(t_2)-F(t_1))^4}
 \end{eqnarray}
  where $G(x)$ is the distribution function corresponding to the density function of $$g(x)=\exp\big(3\eta(\theta)\tau(x)+3Q(\theta)+3R(x)\big).$$ 
\end{proof}

\section{\bf Some Estimators for the Interval Varextropy}\label{section3}
In this section, we propose some estimators for the interval varextropy defined in \eqref{Khas} and  investigate their some asymptotic properties.
\par $\bullet$  \textbf{ First Estimator}
Let $X_{1:n} \leq X_{2:n}\leq \cdots\leq X_{n:n}$ be the ordered sample and $m$ be an integer satisfying $1\leq m\leq n$. 

Set 
 \begin{eqnarray} \label{7}\nonumber
 IV_{mn}(t_1,t_2)&=& \frac{1}{4(F_n(t_2)-F_n(t_1))^3}\frac{1}{2(n-m)}\sum_{j=1}^{n-m}\bigg(\frac{m}{n+1}(X_{j+m:n}-X_{j:n})^{-1}\bigg)^2 \\ \nonumber &&\times \bigg(I(t_1\leq X_{j:n}\leq t_2)+I(t_1\leq X_{j+m:n}\leq t_2)\bigg)-\frac{1}{4(F_n(t_2)-F_n(t_1))^4}\\ \nonumber &&\times  \bigg[\frac{1}{2(n-m)}\sum_{j=1}^{n-m}(\frac{m}{n+1}\big(X_{j+m:n}-X_{j:n})^{-1}\big)\bigg(I(t_1\leq X_{j:n}\leq t_2) \\  &&+ I(t_1\leq X_{j+m:n}\leq t_2)\bigg)\bigg]^2,
 \end{eqnarray}
 where $F_n(x)=\frac{1}{n}\sum_{i=1}^{n}I(X_i \leq x)$ is the empirical distribution function. In the following theorem, we establish the a.s. consistency of  $IV_{mn}(t_1, t_2)$.
 \begin{theorem}\label{Theo6}
 For all bounded densities $f$ provided $m,n\rightarrow \infty$, $\frac{m}{\log(n)}\rightarrow \infty$ and $\frac{m}{n}\rightarrow 0$, we  have for $0<t_1<t_2<\infty$, $$|IV_{mn}(t_1,t_2)-IV(t_1,t_2)|\rightarrow 0  \  \  a.s.$$
\end{theorem}
\begin{proof}
 Set
  \begin{equation} \label{8} 
T_{mn}= \frac{1}{2(n-m)}\sum_{j=1}^{n-m}\bigg[\frac{m}{n+1}\big(X_{j+m:n}-X_{j:n})^{-1}\bigg]^2\bigg(I(t_1 \leq X_{j:n} \leq t_2) + I(t_1\leq X_{j+m:n}\leq t_2)\bigg),
 \end{equation}
and 
  \begin{equation} \label{9} 
T{^\prime}_{mn}= \bigg[\frac{1}{2(n-m)}\sum_{j=1}^{n-m}\bigg(\frac{m}{n+1}\big(X_{j+m:n}-X_{j:n}\big)^{-1}\bigg) \bigg(I(t_1\leq X_{j:n}\leq t_2)+ I(t_1\leq X_{j+m:n}\leq t_2)\bigg)\bigg]^2.
 \end{equation}
Then, we can write
\begin{equation} \label{10} 
|IV_{mn}(t_1,t_2)-IV(t_1,t_2)| \leq  |I_n|+|II_n|,
 \end{equation}
where
\begin{equation} \label{11} 
I_n=\frac{1}{4(F_n(t_2)-F_n(t_1))^3} T_{mn}-\frac{1}{4(F(t_2)-F(t_1))^3}\int_{t_1}^{t_2}(f(x))^3 \,d{x},
 \end{equation}
and
\begin{equation} \label{12} 
II_n=\frac{1}{4(F_n(t_2)-F_n(t_1))^4} T{^\prime}_{mn}-\frac{1}{4(F(t_2)-F(t_1))^4}\big(\int_{t_1}^{t_2}(f(x))^2 \,d{x}\big)^2.
 \end{equation}
On the other hand,

\begin{equation} \label{13} 
I_n=\frac{1}{4}T_{mn} A_n+\frac{1}{4(F(t_2)-F(t_1))^3} \big(T_{mn}-\int_{t_1}^{t_2}f^3(x) \,d{x}\big),
 \end{equation}
where
\begin{equation} \label{14} 
A_n=\frac{1}{(F_n(t_2)-F_n(t_1))^3}-\frac{1}{(F(t_2)-F(t_1))^3}.
 \end{equation}
But,
\begin{equation} \label{15} 
|A_n| \leq \frac{6 \mathop {\sup} \limits_{x } {|F(x)-F_n(x)|}}{(F_n(t_2)-F_n(t_1))^2(F(t_2)-F(t_1))^2}.
 \end{equation}
The Glivenko$-$Cantelli theorem ensures that 
\begin{equation} \label{16} 
\mathop {\lim} \limits_{n \to \infty }{|A_n|}=0, a.s.
 \end{equation}
 \eqref{13},  \eqref{16}  and Theorem \ref{Theo1} of  \cite{VanES1992} ensure that:
 \begin{equation} \label{17} 
\mathop {\lim} \limits_{n \to \infty }{|I_n|}=0, a.s.
 \end{equation}
Also, 
 \begin{equation} \label{18} 
II_n=\frac{1}{4}T^{\prime}_{mn} B_n+\frac{1}{4(F(t_2)-F(t_1))^4} \bigg(T^{\prime}_{mn}-\big(\int_{t_1}^{t_2}f^2(x) \,d{x}\big)^2\bigg),
 \end{equation}
where
\begin{equation} \label{19} 
B_n=\frac{1}{(F_n(t_2)-F_n(t_1))^4}- \frac{1}{(F(t_2)-F(t_1))^4}.
 \end{equation}
The Glivenko$-$Cantelli theorem ensures that 
\begin{equation} \label{20} 
\mathop {\lim} \limits_{n \to \infty }{|B_n|}=0, a.s.
 \end{equation}
From \eqref{18}, \eqref{20} and Theorem \ref{Theo1} of  \cite{VanES1992}, we can see that 
\begin{equation} \label{21} 
\mathop {\lim} \limits_{n \to \infty }{|II_n|}=0, a.s.
 \end{equation}
\eqref{10}, \eqref{17} and \eqref{21} complete the proof.
\end{proof}

\par $\bullet$  \textbf{Second Estimator}
\rm
Set 
 \begin{equation} \label{22} 
\widetilde{IV}(t_1,t_2)= \frac{1}{4(F_n(t_2)-F_n(t_1))^3}\int_{t_1}^{t_2}{f_n}^3(x)\,d{x}-\frac{1}{4(F_n(t_2)-F_n(t_1))^4}\bigg(\int_{t_1}^{t_2}{f_n}^2(x)\,d{x}\bigg)^2, 
 \end{equation}
 where  $f_n(x)=\frac{1}{nh_n}\sum_{j=1}^{n} K(\frac{x-X_i}{h_n})$ is the kernel density estimation.
 
 In the following theorem, we prove a.s. consistency of $\widetilde{IV}(t_1,t_2)$.
  \begin{theorem}\label{Theo7}
Assume that density function  $f$  have a support on $(0, \infty)$ and $f{^\prime}$ is bounded. Let $K$ be a kernel density function and 
$\int_{-\infty}^{\infty} |t|K(t)\,d{t}<\infty$. 
 Then under conditions   $h_n \rightarrow  0$, $n h_n \rightarrow  \infty$ and $h_n^{-1} n^{-\frac{1}{2}}\big(\log\log n\big)^{\frac{1}{2}} \rightarrow 0$ as $n \rightarrow \infty$,   we can see that for every $0<t_1<t_2< \infty$:
 \begin{equation*}  
\mathop {\lim} \limits_{n \to \infty }{|\widetilde{IV}(t_1,t_2)-IV(t_1,t_2)|}=0,\  a.s. 
 \end{equation*}
\end{theorem}
\begin{proof}
 \begin{equation} \label{23} 
|\widetilde{IV}(t_1,t_2)-IV(t_1,t_2)| \leq |I_n|+|II_n|, 
 \end{equation}
 where
 \begin{equation} \label{24} 
I_n=\frac{1}{4}(\int_{t_1}^{t_2}f_n^3(x) \,d{x}) A_n+\frac{1}{4(F(t_2)-F(t_1))^3} \big(\int_{t_1}^{t_2}f_n^3(x) \,d{x}-\int_{t_1}^{t_2}f^3(x) \,d{x}\big),
 \end{equation}
and
 \begin{equation} \label{25} 
II_n=\frac{1}{4}\big(\int_{t_1}^{t_2}f_n^2(x) \,d{x}\big) ^2 B_n+\frac{1}{4(F(t_2)-F(t_1))^4}\bigg(\big[\int_{t_1}^{t_2}f_n^2(x) \,d{x}\big]^2-\big[\int_{t_1}^{t_2}f^2(x) \,d{x}\big]^2\bigg).
 \end{equation}
We use the following inequality in the proof of the theorem.

Let $1 \leq p < \infty$, then for functions $q$ and $u$ in $L_p$, we have
 \begin{eqnarray} \label{26}\nonumber
 \int_{0}^{\infty}\left| |q(t)|^p-|u(t)|^p\right| \,d{\mu(t)} &\leq &  p2^{p-1}  \int_{0}^{\infty}|q(t)-u(t)|^p\,d{\mu(t)}+p2^{p-1} \bigg( \int_{0}^{\infty}|u(t)|^p\,d{\mu(t)}\bigg)^{1-\frac{1}{p}}\\  &\times & \bigg( \int_{0}^{\infty}|q(t)-u(t)|^p\,d{\mu(t)}\bigg)^{\frac{1}{p}},
 \end{eqnarray}
 where $\mu$ is a measure on the Borel sets on $\Re$. Using this inequality we observe that
  \begin{eqnarray} \label{27}\nonumber
\left | \int_{t_1}^{t_2}f_n^3(x) \,d{x}-\int_{t_1}^{t_2}f^3(x) \,d{x}\right | &\leq &  \int_{0}^{\infty} \left | f_n^3(x) -f^3(x) \right |\,d{x}  \leq    12 \int_{0}^{\infty} \left | f_n(x) -f(x)  \right |^3\,d{x} \\ &+& 12 \bigg(\int_{0}^{\infty}f^3(x)\,d{x} \bigg)^{\frac{2}{3}}\bigg(\int_{0}^{\infty}\left|f_n(x)-f(x)\right|^3\,d{x} \bigg)^{\frac{1}{3}}.
 \end{eqnarray}
However,
  \begin{eqnarray} \label{28}\nonumber
 \int_{0}^{\infty} \left | f_n(x) \,d{x}-f(x) \right |^3\,d{x}  &\leq&  \big(\mathop {\sup} \limits_{x } |f_n(x)-f(x)|^2\big) \int_{\Re} \left | f_n(x)- f(x) \right| \,d{x}\\&:=&Y_{n,1}  \times Y_{n,2}.
 \end{eqnarray}
Now, by using the Chung's law of the iterated Logarithms (LIL), we can write
  \begin{eqnarray} \label{29}\nonumber
\mathop {\sup} \limits_{x } |f_n(x)-f(x)|&\leq & \mathop {\sup} \limits_{x } |f_n(x)-E(f_n(x))|+\mathop {\sup} \limits_{x } |E(f_n(x))-f(x)| \\ \nonumber &\leq&
\mathop {\sup} \limits_{x } \left |h_n^{-1} \int_{\Re} K(\frac{x-y}{h_n})\,d{F_n(y)}-h_n^{-1}\int_{\Re} K(\frac{x-y}{h_n})\,d{F(y)} \right |\\ \nonumber &+&\mathop {\sup} \limits_{x } \left |\int_{\Re} \big(f(x-h_n t)-f(x)\big)K(t)\,d{t}\right| \leq  \mathop {\sup} \limits_{x } \left |f^{\prime}(x) \right| h_n^{-1}  \mathop {\sup} \limits_{x } |F_n(x)-F(x)|\\ \nonumber &+& \mathop {\sup} \limits_{x } \left |f^{\prime}(x) \right| h_n  \int_{\Re} |t|K(t)\,d{t}\\ &=&h_n^{-1}O\big(n^{-\frac{1}{2}}(\log \log n )^{\frac{1}{2}}\big)+O(h_n), a.s.
 \end{eqnarray}
Clearly, by \eqref{29},
 \begin{equation} \label{30} 
Y_{n,1}=O\bigg(\bigg[\bigg(h_n^{-1}n^{-\frac{1}{2}}(\log \log n  )^{\frac{1}{2}}\bigg)\vee h_n\bigg]^2\bigg), a.s.
 \end{equation}
Also, by the results on the a.s. behavior of $L_1$-norms of kernel density estimators (see page 149 of  \cite{EggermontandLariccia2001}),
 \begin{equation} \label{31} 
Y_{n,2}=O\bigg(h_n^{2}+(n h_n)^{-\frac{1}{2}}+(n^{-1}\log(n))^{\frac{1}{2}}\bigg), a.s.
 \end{equation}
\eqref{16},\eqref{24},\eqref{27},\eqref{28},\eqref{30} and \eqref{31}, ensure that 
 \begin{equation} \label{32} 
\mathop {\lim} \limits_{n \to \infty }{|I_n|}=0,\  a.s. 
 \end{equation}
 In a similar way to the above discussion and using formula \eqref{20}, it  can be shown that 
 \begin{equation} \label{33} 
\mathop {\lim} \limits_{n \to \infty }{|II_n|}=0,\  a.s.
 \end{equation}
 \eqref{23},\eqref{32} and \eqref{33} complete the proof.
\end{proof}

\par $\bullet$  \textbf{ Third Estimator}

Set 
 \begin{equation*} 
\widehat{IV}(t_1,t_2)= \frac{1}{4(F_n(t_2)-F_n(t_1))^3}\int_{t_1}^{t_2}{f_n}^2(x)\,d{F_n(x)}-\frac{1}{4(F_n(t_2)-F_n(t_1))^4}\bigg(\int_{t_1}^{t_2}f_n(x)\,d{F_n(x)}\bigg)^2.
 \end{equation*}
 
In the following theorem we show  a.s.  consistency of $\widehat{IV}(t_1,t_2)$ under some conditions.
  \begin{theorem}\label{Theo8}
Let $f$ and $f^{\prime}$ be bounded. If $h_n\rightarrow 0$  as $n\rightarrow \infty$, then for every $0<t_1<t_2<\infty$, we can see that 
  \begin{equation*}
\mathop {\lim} \limits_{n \to \infty }{|\widehat{IV}(t_1,t_2)-IV(t_1,t_2)|}=0,\  a.s. 
 \end{equation*}
\end{theorem}

\begin{proof}
  \begin{equation}\label{34}
|\widehat{IV}(t_1,t_2)-IV(t_1,t_2)| \leq |I_n|+|II_n|,
 \end{equation}
 where
  \begin{equation}\label{35}
I_n=\frac{1}{4} A_n \int_{t_1}^{t_2}f_n^2  \,d{F_n(x)} +\frac{1}{4(F(t_2)-F(t_1))^3} \bigg[\int_{t_1}^{t_2}f_n^2  \,d{F_n(x)}-\int_{t_1}^{t_2}f^2(x) \,d{F(x)}\bigg],
 \end{equation}
 and
   \begin{equation}\label{36}
II_n=\frac{B_n}{4} \big( \int_{t_1}^{t_2}f_n \,d{F_n(x)}\big)^2 -\frac{1}{4(F(t_2)-F(t_1))^4} \bigg(\bigg[\int_{t_1}^{t_2}f(x) \,d{F(x)}\bigg]^2-\bigg[\int_{t_1}^{t_2}f_n(x) \,d{F_n(x)}\bigg]^2\bigg).
 \end{equation}
 On the other hand,
   \begin{eqnarray} \label{37}\nonumber
\left |\int_{t_1}^{t_2}f_n^{2} \,d{F_n(x)}-\int_{t_1}^{t_2}f^2(x) \,d{F(x)}\right| &\leq &\int_{t_1}^{t_2}\left |f_n^{2}(x)-f^2(x)\right| \,d{F_n(x)}+\left| \int_{t_1}^{t_2}f^2(x) \,d\big({F_n(x)}-F(x)\big)\right| \\ \nonumber &\leq & \mathop {\sup} \limits_{x } \left |f_n^{2}(x)-f^2(x) \right|(F_n(t_2)-F_n(t_1))\\ \nonumber &+&  \left| \int_{t_1}^{t_2}\big(F_n(x)-F(x)\big)f(x)f^{\prime}(x) \,d{x}\right| \\ \nonumber & \leq & \mathop {\sup} \limits_{x } \left |f_n(x)-f(x) \right|^2+2(\mathop {\sup} \limits_{x } f(x)) \mathop {\sup} \limits_{x }\left |f_n(x)-f(x) \right|\\  & + & \mathop {\sup} \limits_{x }\left |F_n(x)-F(x) \right|(\mathop {\sup} \limits_{x } f(x))\mathop {\sup} \limits_{x } |f^{\prime}(x)| (t_2-t_1),
 \end{eqnarray}
the Chung's LIL, \eqref{37} and \eqref{29} ensure that 
\begin{equation}\label{38}
\mathop {\lim} \limits_{n \to \infty }{\left |\int_{t_1}^{t_2}f_n^{2} \,d{F_n(x)}-\int_{t_1}^{t_2}f^2(x) \,d{F(x)}\right|}=0, a.s.
\end{equation}
From \eqref{16}, \eqref{35} and \eqref{38}, we have 
  \begin{equation}\label{39}
\mathop {\lim} \limits_{n \to \infty }{|I_n|}=0, a.s.
 \end{equation}
 In a similar way we can see that
 \begin{equation}\label{40}
\mathop {\lim} \limits_{n \to \infty }{|II_n|}=0, a.s.
 \end{equation}
 \eqref{34}, \eqref{39} and \eqref{40} complete the proof.
\end{proof}

\rm
\section{\bf Simulation Study}\label{section4}
In this section, the results of the Monte Carlo studies on biases and mean squared errors (MSEs) of our proposed estimators are presented. We consider the gamma, uniform and exponential distributions, which are considered in many references. For each sample size, 10,000 samples were generated and MSEs  and biases of estimators were computed. We used the $m=[\sqrt n+0.5]$  formula to estimate the  interval varextropy suggested by \cite{GrzegorzewskiandWieczorkowski1999}. They used this heuristic formula for entropy estimation. In addition, we choose the  Epanechnikov kernel $K(x)=\frac{3}{4}(1-x^2)I(|x|<1)$  as the kernel and its corresponding optimal value of $h_n=1.06sn^{-\frac{1}{5}}$, where $s$ is the sample standard deviation. Tables \ref{Tab1}-\ref{Tab3} give simulated biases and MSEs  of $IV_{m, n}(t_1, t_2)$, $\widetilde{IV}(t_1,t_2)$ and $\widehat{IV}(t_1,t_2)$  for gamma, uniform and exponential distributed samples, respectively. 
For gamma and  exponential distributions, we use $t_1=0$ and $t_2=3$ and for uniform distribution, we use $t_1=0$ and $t_2=0.5$.
\begin{table}[h]
\caption{MSE and Bias of estimators for estimating  the  interval varextropy $IV(0,3)$ of the gamma distribution ($G(2,1)$).}\label{Tab1}
\centering
\begin{tabular}{lccc}\hline
\multicolumn{4}{c}{MSE(Bias)}\\
\hline
$n$ &{$IV_{m, n}(0,3)$} &{$\widetilde{IV}(0,3)$} &{$\widehat{IV}(0,3)$}\\
\hline
10 &   0.4952(0.4952) &{\bf 0.2426(-0.2277)}&0.3505(-0.3491)  \\
20 & 0.4551(0.4530)  & {\bf 0.1464(-0.0684) }& 0.2239(-0.2163)\\
30 & 0.2067(0.0667) &{\bf 0.1300(-0.0533)}& 0.1615(-0.1228)  \\
40&0.0405(0.0398)&{\bf 0.0384(-0.0342)}& 0.0457(-0.0343) \\
50&0.0343(0.0265) &{\bf 0.0324(-0.0254)} &0.0406(-0.0335) \\ 
100 & 0.0298(0.0213) & {\bf 0.0286(-0.0150)} & 0.0290(-0.0251)\\
\hline
\end{tabular}
\end{table}

\begin{table}[h]
\caption{MSE and Bias of estimators for estimating  the  interval varextropy $IV(0,0.5)$ of the uniform distribution ($U(0,1)$).}\label{Tab2}
\centering
\begin{tabular}{lccc}\hline
\multicolumn{4}{c}{MSE(Bias)}\\
\hline
$n$ &{$IV_{m, n}(0,0.5)$} &{$\widetilde{IV}(0,0.5)$} &{$\widehat{IV}(0,0.5)$}\\
\hline
10 &   0.4685(0.4679) & 0.2278(0.2272)&{\bf 0.0242(0.0227)} \\
20 & 0.3314(0.3894)  & 0.1263(0.0985) &{\bf  0.0146(0.0048)}\\
30 & 0.3140(0.3100) & 0.1093(0.0906)&{\bf 0.0130(0.0033)}  \\
40&0.0966(0.0264)& 0.0475(0.0149)&{\bf 0.0038(0.0031)} \\
50&0.0873(0.0216) & 0.0433(0.0117) &{\bf 0.0031(0.0030)} \\ 
100 & 0.0815(0.0120) &  0.0406(0.0066) &{\bf 0.0021(0.0027)}\\
\hline
\end{tabular}
\end{table}

\begin{table}[h]
\caption{MSE and Bias of estimators for estimating  the  interval varextropy $IV(0,3)$ of the exponential distribution ($exp(1)$).}\label{Tab3}
\centering
\begin{tabular}{lccc}\hline
\multicolumn{4}{c}{MSE(Bias)}\\
\hline
$n$ &{$IV_{m, n}(0,3)$} &{$\widetilde{IV}(0,3)$} &{$\widehat{IV}(0,3)$}\\
\hline
10 &   0.4557(0.4544) &{\bf 0.2321(-0.2250)}&0.3417(-0.3374)  \\
20 & 0.3526(0.3571)  & {\bf 0.1002(-0.0193) }& 0.2085(-0.2080)\\
30 & 0.3187(0.3413) &{\bf 0.0949(-0.0187)}& 0.2023(-0.1913)  \\
40&0.0951(0.0298)&{\bf 0.0623(-0.0174)}& 0.0626(-0.0181) \\
50&0.0866(0.0246) &{\bf 0.0567(-0.0138)} &0.0569(-0.0141) \\ 
100 & 0.0813(0.0193) & {\bf 0.0503(-0.0122)} & 0.0517(-0.0130)\\
\hline
\end{tabular}
\end{table}

It is apparent from the tables that  both the absolute bias and MSE  decrease as the sample size increases. The bold type in these tables indicates the  interval varextropy estimator achieves the minimal MSE (bias). In the case of gamma distribution, Table \ref{Tab1}, it is observed that, the estimator $\widetilde{IV}(0,3)$ has the best performance in terms of bias and MSE. Under the uniform distribution, Table \ref{Tab2}, it can be seen that, the MSEs and biases of $\widehat{IV}(0,0.5)$ are always smaller than other proposed estimators. 

According to Table \ref{Tab3},  under exponential distribution, the estimator $\widetilde{IV}(0,3)$ performs better than other proposed estimators in terms of MSE and bias. 
\section{\bf Some Tests of Uniformity}\label{section5}
In this section, we introduce some goodness-of-fit tests of uniformity. These tests are based on our proposed interval varextropy estimators and their percentage points and power values are obtained by Monte Carlo simulation.

Consider the class of continuous distribution functions ${F}$ with density function $f(x)$ concentrated on the interval $[0, 1]$. 
Since the interval varextropy of a random variable is non-negative, we have $IV(t_1,t_2)=var(-\frac{1}{2}f(X_{t_1,t_2}))\geq 0$.

If  $f$ is a standard uniform density function, it is easy to see that for each  $(t_1, t_2) \in D=\{(u,v)\in  \Re _{+}^{2}: F(u)<F(v)\}$, $IV(t_1, t_2)=0$.

On the other hand, if there exists a $(t_1, t_2) \in D$ where $IV(t_1, t_2)=0$, then we have $f(x_{t_1, t_2})=\frac{f(x)}{F(t_2)-F(t_1)}=c$ on $[0, 1]$, where $c$ is a constant.  From the property of density function $f(x_{t_1, t_2})$, it is obvious that $f(x_{t_1, t_2})=c=t_2-t_1$, where implies that 
\begin{equation}\label{41}
f(x)=(t_2-t_1)(F(t_2)-F(t_1)),   \ \forall   x \in [0,1].
 \end{equation}
Also, the Mean Value theorem ensures that there exists  a $l \in (t_1, t_2)$, such that 
\begin{equation}\label{42}
f(l)=\frac{F(t_2)-F(t_1)}{t_2-t_1}.
 \end{equation}
 \eqref{41} and \eqref{42} imply that $f(l)=(t_2-t_1)(F(t_2)-F(_1))=\frac{F(t_2)-F(t_1)}{t_2-t_1}$, i.e.,  $t_2-t_1=1$ or $t_2=1$, $t_1=0$. 
 
 This result and \eqref{41} ensure that $f$ is a standard uniform density function.
  
The above discussion gives us the idea to use the proposed interval  varextropy estimators as the test statistics of goodness-of-fit tests of uniformity. 

Let $X_1$, $\ldots$, $X_n$ be a random sample from a continuous distribution function $F(x)$ on $[0,1]$. The null hypothesis is $H_0: F(x)=x$ and the alternative hypothesis (denoted $H_1$) is the opposite of $H_0$. Given any significance level $\alpha$, our hypothesis-testing procedure can be defined by the critical region:  
\begin{equation}\label{43}
G_n=IV_n(t_1, t_2)\geq C_{1-\alpha}, 
\end{equation}
where $IV_n(t_1, t_2)$ is one of the proposed estimators and $C_{1-\alpha}$ is the critical value for the test with level $\alpha$. Notice when $IV_n(t_1, t_2)\stackrel{p}{\to}IV(t_1, t_2)$, under $H_0$: $G_n\stackrel{p}{\to}0$ and $H_1$: $G_n$ converges  to a number larger than zero in probability.

Due to the complexity of the calculations, it is not easy to determine the distribution of the test statistics under $H_0$, so the critical values are calculated using the Monte Carlo method. Based on three different interval  varextropy estimators, we can propose the following test statistics for testing the uniformity:  
\begin{eqnarray}\nonumber
GV_{n}&=IV_{m, n}(0, 1),\\ \nonumber
GD_n&=\widetilde{IV}(0,1),\\ \nonumber
GB_n&=\widehat{IV}(0,1).
\end{eqnarray}
We compare the powers of our proposed test statistics with Kolmogorov-Smirnov statistic (see \cite{Kolmogorov1933} and \cite{Smirnov1939}).
The test statistic of Kolmogorov-Smirnov is given as follows:
\begin{equation}\label{44}
KS=max\left(max_{1\leq i\leq n}\Big\{\frac{i}{n}-X_{(i)}\Big\}, max_{1\leq i\leq n}\Big\{X_{(i)}-\frac{i-1}{n}\Big\}\right),
\end{equation}
where $X_{(1)}$,  $\ldots$, $X_{(n)}$ are order statistics. 

\begin{table}[ht] 
 \small
 \caption{Percentage points of the proposed test statistic at the level $\alpha=0.05.$}\label{Tab4}
\begin{center}
 \begin{tabular}{cccccc}
 \hline
 n&$GV_{n}$&$GD_n$&$GB_n$\\
 \hline
 10 & 0.9102 & 0.0665&0.0574\\
 20 &0.4937&0.0485&0.0374\\
30&0.2845&0.0413&0.0297\\
40&0.1933&0.0371&0.0256\\
50&0.1389&0.0343&0.0226\\
75&0.0872&0.0296&0.0183\\
100&0.0715&0.0270&0.0158\\
\hline
\end{tabular}
\end{center}
\end{table}

\begin{table}[ht] 
\small
\caption{Power comparisons of the proposed tests at a significance level of 0.05. }\label{Tab5}
\begin{center}
\begin{tabular}{|c|c|c|c|c|c|}
\hline
 n& $Alternative$&$GV_{n}$&$GD_n$& $GB_n$&$KS$\\
 \hline
&$A_{1.5}$ & 0.0739&\bf 0.1767&0.1072& 0.1542\\
&$A_2$ & 0.1150&\bf 0.3959&0.2033& 0.3871 \\
&$B_{1.5}$&0.0719&0.1589&{\bf 0.1629}&0.0367\\
10&$B_2$&0.1126&{\bf 0.3494}&0.3292&0.0447\\
&$B_3$&0.2156&{\bf 0.7326}&0.6421&0.0873\\
&$C_{1.5}$& 0.0808&0.0415&0.0296&{\bf 0.1148}\\
&$C_2$& 0.1295&0.0406&0.0300&{\bf 0.2001}\\
\hline
\hline
&$A_{1.5}$&0.0884&{\bf 0.2889}&0.1575&0.2796\\
&$A_2$&0.1684& 0.6680&0.3648&{\bf 0.6979}\\
&$B1.5$&0.0879&0.2424&{\bf 0.2543}& 0.0579\\
20&$B_2$&0.1640&{\bf 0.5707}&0.5493&0.1246\\
&$B_3$&0.3689&{\bf 0.9459}&0.8948&0.4133\\
&$C_{1.5}$&0.09628&0.0373&0.0251& {\bf 0.1501}\\
&$C_2$& 0.1900&0.0440&0.0281&{\bf 0.3159}\\
\hline
\hline
&$A_{1.5}$&0.1114& 0.3933&0.2092&{\bf 0.4066}\\
&$A_2$&0.2425&0.8319&0.5078&{\bf 0.8751}\\
&$B1.5$&0.1100&0.3160&{\bf 0.3416}& 0.0837\\
30&$B_2$&0.2333&{\bf 0.7251}&0.7132&0.2411\\
&$B_3$&0.5500&{\bf 0.9907}&0.9736&0.7524\\
&$C_{1.5}$&0.1236&0.0385&0.0228&{\bf  0.1892}\\
&$C_2$& 0.2827&0.0498&0.0263&{\bf 0.4383}\\
\hline 
\end{tabular}
\end{center}
\end{table}
For power comparisons, we compute the powers of  our proposed tests and the powers of tests mentioned above under the following alternative distributions:
\begin{eqnarray*}
A_k&: F(x)=1-{(1-x)}^k, \hspace{0.5cm}0\leq x\leq 1\mbox\ \ {({\rm for}\ k=1.5, 2)};\\
B_k&: F(x)=\left\{ 
\begin{array}{llc}
2^{k-1}x^k,&0\leq x\leq 0.5&\\
&&\mbox\ \ {({\rm for}\ k=1.5, 2, 3)};\\
1-2^{k-1}{(1-x)}^k,& 0.5\leq x\leq 1&
\end{array}\right.
\\
C_k&: F(x)=\left\{ 
\begin{array}{llc}
0.5-2^{k-1}{(0.5-x)}^k,&0\leq x\leq 0.5&\\
&&\mbox\ \ {({\rm for}\ k=1.5, 2).}\\
0.5+2^{k-1}{(x-0.5)}^k,& 0.5\leq x\leq 1&
\end{array}\right.
\end{eqnarray*}
\cite{Stephens1974} uses these alternative distributions in his study on power comparisons of some uniformity tests. 

Based on $100,000$ repetitions, the critical values are estimated as shown in Table \ref{Tab4}.
Table \ref{Tab5} contains the results of simulations of our proposed test statistics and competitive test statistic with $100,000$ repetitions for $n=10, 20, 30$ and $\alpha=0.05$. The bold type in Table \ref{Tab5} indicates the statistic achieving the maximal power. 
According to Table \ref{Tab5}, the performance  of tests depends on alternative distributions.
 
Against alternative $A$, our proposed test $GD_{n}$ and  $KS$-based test  give competitive powers to different samples.  For sample size $n=10$, the test based on   $GD_{n}$ is the best. For moderate sample size $n=20$, if alternative distribution is $A_{1.5}$, $GD_{n}$ is better. For 
$n=30$, $KS$ test   performs better than the others. For alternative distribution $B$, $GD_{n}$ and $GB_{n}$ are the best tests and other tests have low powers in this case.  If the alternative is $C$, the test based on $KS$ is the best test.
\section{\bf Real Data}\label{section6}
In this section, we illustrate the use of proposed estimators in the real case. The data set considered here represents the remission times (in months) of a random sample of 128 bladder cancer patients. These data were analyzed by \cite{LeeandWang2003}.
The data are:

0.08, 2.09, 3.48, 4.87, 6.94, 8.66, 13.11, 23.63, 0.20, 2.23, 3.52, 4.98, 6.97, 9.02, 13.29,
0.40, 2.26, 3.57, 5.06, 7.09, 9.22, 13.80, 25.74, 0.50, 2.46, 3.64, 5.09, 7.26, 9.47, 14.24, 25.82, 0.51,
2.54, 3.70, 5.17, 7.28, 9.74, 14.76, 6.31, 0.81, 2.62, 3.82, 5.32, 7.32, 10.06, 14.77, 32.15, 2.64, 3.88,
5.32, 7.39, 10.34, 14.83, 34.26, 0.90, 2.69, 4.18, 5.34, 7.59, 10.66, 15.96, 36.66, 1.05, 2.69, 4.23, 5.41,
7.62, 10.75, 16.62, 43.01, 1.19, 2.75, 4.26, 5.41, 7.63, 17.12, 46.12, 1.26, 2.83, 4.33, 5.49, 7.66, 11.25,
17.14, 79.05, 1.35, 2.87, 5.62, 7.87, 11.64, 17.36, 1.40, 3.02, 4.34, 5.71, 7.93, 11.79, 18.10, 1.46, 4.40,
5.85, 8.26, 11.98, 19.13, 1.76, 3.25, 4.50, 6.25, 8.37, 12.02, 2.02, 3.31, 4.51, 6.54, 8.53, 12.03, 20.28,
2.02, 3.36, 6.76, 12.07, 21.73, 2.07, 3.36, 6.93, 8.65, 12.63, 22.69.

\cite{Shankeretal2015} fitted this data using an exponential distribution with $\lambda=0.106773$. Considering this model, we have obtained the values of $IV(t_1,t_2)$ and estimators for different $(t_1,t_2)$ in Table \ref{Tab6}. The closeness of the estimated results to the values of  $IV(t_1,t_2)$ shows that our estimators are good in the real scenario.
\begin{table}[h]
\caption{$IV(t_1,t_2)$ and their estimators for cancer data under different truncation limits.}\label{Tab6}
\centering
\begin{tabular}{lcccc}\hline

\hline
$(t_1,t_2)$ &{$IV_{m, n}(t_1,t_2)$} &{$\widetilde{IV}(t_1,t_2)$}&{$\widehat{IV}(t_1,t_2)$} &{$IV(t_1,t_2)$}\\
\hline
(1,7) &   0.01153032 &0.0006164359&0.0007132116&0.0003462006  \\ \hline
(1,13) & 0.00316791  & 0.0002297752& 0.0005251652&0.0002789512\\ \hline
(2,10)& 0.00335218 & 0.0004016302& 0.0006531906&0.0002586823\\

\hline
\end{tabular}
\end{table}

\section{\bf Concluding Remarks}\label{section7}
In this paper, we introduced the concept of interval varextropy where  is a varextropy measure for doubly truncated random variables. We investigated changes of this measure under linear transforms.
We showed that if $IV(t_1,t_2)$ is constant, then $IJ(t_1,t_2)$ can be obtained in terms of GFR functions. Also, we showed that if $IV(t_1,t_2)$ is constant, then under some conditions, the desired distribution is exponential. We obtained lower and upper bounds for $IV(t_1,t_2)$.
We proposed some estimators for $IV(t_1,t_2)$ and proved some asymptotic properties for these estimators.

In the simulation study, we investigated the behavior of estimators based on bias and MSE. By using the proposed estimators, we proposed some tests for uniformity tests. Finally, we evaluated estimators for a  real data set.
 \section*{Acknowledgments}
The authors would like to thank the referees for their careful reading and their useful comments which led to this considerably improved version.
{ }
\end{document}